\def\Bbb#1{{\bf #1}}

\def\fnote#1{\footnote}
\def\blacksquare{\hbox{\vrule width 4pt height 4pt depth 0pt}}

\def\cwleftpar#1#2{\leftskip #1 \rightskip #2 plus 1fill}
\def\cwrightpar#1#2{\leftskip #1 plus 1fill \rightskip #2}
\def\cwcenterpar#1#2{\leftskip #1 plus 1fill \rightskip #2 plus 1fill}
\def\cwfullpar#1#2{\leftskip#1\rightskip#2}

\def\cwoutdent#1#2{\llap{\hbox to #1{#2 \hss}}\ignorespaces}
\def\cwparbegin#1#2#3#4#5{
	\ifcase #1 \cwleftpar{#2}{#3}
	\or \cwrightpar{#2}{#3}
	\or \cwcenterpar{#2}{#3}
	\else \cwfullpar{#2}{#3}\fi
	\ifcase #4 \baselineskip = 1.5\baselineskip
	\or \baselineskip = 2\baselineskip
	\or \baselineskip = 3\baselineskip
	\else \baselineskip = 1\baselineskip\fi
	\ifdim #5 > 0in \else \noindent \fi
	\noindent\ignorespaces}
\documentclass{article}
\begin{document}
\advance \vsize by -1\baselineskip
\def\makefootline{
{\vskip \baselineskip \noindent \folio                                  \par
}}

\vspace*{2ex}

\noindent {\Huge Linear Transports along Paths in\\[0.4ex]
 					Vector Bundles}\\[1.8ex]
\noindent {\Large IV. Consistency with Bundle Metrics}
\vspace*{2ex}

\noindent Bozhidar Zakhariev Iliev
\fnote{0}{\noindent $^{\hbox{}}$Permanent address:
Laboratory of Mathematical Modeling in Physics,
Institute for Nuclear Research and \mbox{Nuclear} Energy,
Bulgarian Academy of Sciences,
Boul.\ Tzarigradsko chauss\'ee~72, 1784 Sofia, Bulgaria\\
\indent E-mail address: bozho@inrne.bas.bg\\
\indent URL: http://theo.inrne.bas.bg/$\sim$bozho/}

\vspace*{2ex}

{\bf \noindent Published: Communication JINR, E5-94-17, Dubna, 1994}\\[1ex]
\hphantom{\bf Published: }
http://www.arXiv.org e-Print archive No.~math.DG/0505008\\[2ex]

2003 PACS numbers: 02.40.Ma, 02.40.Vh, 04.90.+e\\[2ex]

\noindent
{\small
The \LaTeXe\ source file of this paper was produced by converting a
ChiWriter 3.16 source file into
ChiWriter 4.0 file and then converting the latter file into a
\LaTeX\ 2.09 source file, which was manually edited for correcting numerous
errors and for improving the appearance of the text.  As a result of this
procedure, some errors in the text may exist.
}\\[1ex]

	\begin{abstract}
 The problem for consistency between linear transports along paths and real
bundle metrics in real vector bundles is stated. Necessary and/or sufficient
conditions, as well as conditions for existence, for such consistency are
derived. All metrics (resp. transports) consistent with a given transport
(resp. metric) are explicitly obtained. The special case of linear transports,
generated by derivations of tensor algebras, of vectors is considered.
Analogous problems are investigated in complex vector bundles endowed with
Hermitian metrics.
	\end{abstract}\vspace{3ex}

{\bf 1. INTRODUCTION}
\medskip

The present work investigates certain connections between linear transports along paths and bundle metrics, defined on one and the same vector bundle, which arise from the question of their consistency. More precisely, we shall consider the following main problem.

Let $(E,\pi ,B)$ be a real vector bundle with a base $B$, total bundle space
$E$ and projection $\pi :E \to B [1]$. Let $g$ be a bundle metric on it [7],
i.e. it is a map $g:x\mapsto g_{x}, x\in B$, where the map
\[
 g_{x}:\pi ^{-1}(x) \pi ^{-1}(x) \to {\Bbb R}\qquad (1.1)
\]
 is bilinear,
nondegenerate and symmetric for every $x\in $B. By definition the scalar
product of $u,v\in \pi ^{-1}(x)$ is $u\cdot v:=g_{x}(u,v), x\in $B. Let $L$
be a linear transport $(L$-transport) along paths in $(E,\pi ,B)$, i.e. if
$\gamma :J \to B, J$ being an arbitrary real interval, is a path in $B$, then
$L:\gamma \mapsto L^{\gamma }$, where $L^{\gamma }$is the $L$-transport along
$\gamma $ and $L^{\gamma }:(s,t)\mapsto L^{\gamma }_{s\to t}, s,t\in
J$ is the $L$-transport along $\gamma $ from $s$ to $t$ having the described
in [2] properties.

{\bf Definition 1.1.} The transport along paths $L$ and the bundle metric
$g$ are consistent (resp. along a path $\gamma :J\to B)$ if $L$ preserves the
scalar products of the vectors along every (resp. along the given) path
$\gamma :J\to B$, i.e. if the scalar product of $u,v\in \pi ^{-1}(\gamma
(s)), s\in J$ is equal to the scalar product of the vectors obtained from $u$
and $v$ by $L$-transportation along $\gamma $ at an arbitrary point $\gamma
(t), t\in J$:
\[
g_{\gamma (s)}(u,v)=g_{\gamma (t)}{\bigl(}L^{\gamma }_{s\to
t}u,L^{\gamma }_{s\to t}v{\bigr)},\quad  s,t\in J.
 \qquad (1.2)
\]
 Due to the arbitrariness of $u$ and $v$ and the nondegeneracy of $g$ the
equality (1.2) is equivalent to
\[
g_{\gamma (s)}=g_{\gamma (t)}\circ {\bigl(}L^{\gamma }_{s\to t}
L^{\gamma }_{s\to t}{\bigr)},\quad s,t\in J.
 \qquad (1.3)
\]
Important examples for transports along paths (in this case along curves)
consistent with some metric are the parallel and Fermi-Walker transports in a
Riemannian manifold which are consistent with the defining them Riemannian
metric. This is proved, for instance, in [3,4] (see also below section 3).

The purpose of this work is to find necessary and/or sufficient conditions for consistency, or one may say compatibility, between $L$-transports along paths and bundle metrics. Analogous problems have been investigated (in a slightly different notation) in [5], where they were studied in the special case of the tangent bundle to a differentiable manifold.

The organization of the material is the following. In Sect. 2 attention is focussed on finding necessary and/or sufficient conditions for local (i.e. along a fixed path) or global (i.e. along any path) consistency between $L$-transports along paths and real bundle metrics. Also the corresponding problems of existence are considered. In Sect. 3 are investigated problemso of consistency concerning the specific case of generated by derivations of tensor algebras $L$-transports of vectors. In Sect. 4  the results of Sect. 2 are transferred to complex vector bundles endowed with Hermitian bundle metrics. Some remarks on the presented in this paper material are made in Sect. 5.

At the end of this introduction we shall note that in a  field of bases
$\{e_{i}, i=1,\ldots  ,\dim(\pi ^{-1}(x)), x\in B\}$ along a  path  $\gamma
:J \to B, J$ being an arbitrary ${\Bbb R}$-interval, any $L$-transport along
$\gamma $ from $s$ to $t, s,t\in J$ is uniquely characterized by its matrix
$H(t,s;\gamma )= H^{i}_{.j}(t,s;\gamma ) $. Hereafter the Latin indices run
from 1 to $n:=\dim(\pi ^{-1}(x)), x\in B$  and further the usual summation
rule from 1 to $n$ is  assumed  over  the repeated on different levels
indices. The general form of  $H(t,s;\gamma )$ is
\[
 H(t,s;\gamma )=F^{-1}(t;\gamma )F(s;\gamma )\qquad (1.4)
\]
 for a
nondegenerate matrix function $F$, defined up to a left multiplication with a
nondegenerate matrix depending only on $\gamma $, i.e. up to the
transformation
\[
F(s;\gamma )\mapsto D(\gamma )F(s;\gamma ), \det(D(\gamma ))\neq 0,\infty
.\qquad (1.5)
\]
For further details of notation and results concerning
$L$-transports along paths in vector bundles the reader is referred to [2].

\medskip
\medskip
 {\bf 2. GENERAL CONDITIONS FOR CONSISTENCY}

\medskip
Let in the real vector bundle $(E,\pi ,B)$ be given an $L$-transport along
paths $L$ and a real bundle metric g. Let $\gamma :J\to B$ and in every fibre
$\pi ^{-1}(\gamma (s)), s\in J$ be fixed a basis $\{e_{i}(s): i=1,\ldots
,n=\dim(\pi ^{-1}(x)), x\in B\}$. Let in $\{e_{i}\}$ the transport
$L^{\gamma }_{s\to t}$along $\gamma $ from $s$ to $t$ be given by its matrix
$H(t,s;\gamma )= H^{\hbox{i.}}_{.j}(t,s;\gamma )  [2]$.

 The components of the metric $g$ in $\{e_{i}(s)\}$ at $\gamma (s)$ are
\[
(g_{\gamma (s)})_{ij}:=g_{\gamma (s)}(e_{i}(s),e_{j}(s)).\qquad (2.1)
\]
 Let
$G(\gamma (s)):=[ (g_{\gamma (s)})_{ij}] $. The nondegeneracy of $g$ means
that $\det(G(\gamma (s)))\neq 0,\infty , s\in J$, and its symmetry means that
$G^{T}(\gamma (s))=G(\gamma (s))$, $s\in J$, where the superscript $T$ means a
transposition of matrices.

{\bf Proposition} ${\bf 2}{\bf .}{\bf 1}{\bf .} A$ necessary and sufficient
condition for a global (resp. local) consistency between a bundle metric $g$
and an $L$-transport along paths $L$ is
\[
G(\gamma (s))
=H^{T}(t,s;\gamma ) G(\gamma (t)) H(t,s;\gamma ),\quad s,t\in J\qquad (2.2)
\]
 for every (resp. for the given) path $\gamma :J\to B$.

{\bf Proof.} If $u=u^{i}e_{i}(s)$ and $v=v^{i}e_{i}(s) (a$ summation from 1
to $n$ is understood over repeated on different levels indices) are vectors
from $\pi ^{-1}(\gamma (s))$, then, due to the bilinearity of the metric, we
have
\[
g_{\gamma (s)}(u,v)=u^{i}v^{j}g_{\gamma
(s)}(e_{i}(s),e_{j}(s))=u^{i}v^{j}(g_{\gamma (s)})_{ij}.\qquad (2.3)
\]
Substituting this equality into (1.2) and taking into account the linearity
of $L$, the definition of its matrix $H [2]$, and the arbitrariness of $u$
and $v$, we see (1.2) to be equivalent to
\[
 (g_{\gamma (s)})_{ij}=H^{\hbox{k.}}_{.i}(t,s;\gamma )(g_{\gamma
(t)})_{kl}H^{\hbox{l.}}_{.j}(t,s;\gamma ),\qquad (2.4)
\]
 which is simply a
component form of (2.2).\blacksquare

Using the general form of the matrices of $L$-transports along paths, which is given by (1.4), we can "simplify" (2.2) by putting it into "one point" form which is the contents of

{\bf Proposition 2.2.} If the nondegenerate matrix function $F(s;\gamma )$
defines the $L$-transport along paths $L$ by (1.4), then a necessary and
sufficient condition for a global (resp. local) consistency of $L$ with the
bundle metric $g$ is the existence of a nondegenerate, symmetric (resp. and
constant) matrix function $C$ of $\gamma $ such that
\[
{\bigl(}F^{-1}(s;\gamma ){\bigr)}^{T}G(\gamma (s))F^{-1}(s;\gamma
)=C(\gamma ),\quad  s\in J,\ C^{T}=C.
 \qquad (2.5)
\]

 {\bf Proof.} Substituting
(1.4) into (2.2) and multiplying the so obtained equality on the left by
${\bigl(}F^{-1}(s;\gamma ){\bigr)}^{T}$and on the right by
$F^{-1}(s;\gamma )$, we find that (2.2) is equivalent to
\[
{\bigl(}F^{-1}(s;\gamma ){\bigr)}^{T}G(\gamma (s))F^{-1}(s;\gamma
)={\bigl(}F^{-1}(t;\gamma ){\bigr)}^{T}G(\gamma (t))F^{-1}(t;\gamma ).
\qquad (2.6)
\]
If the metric and transport are consistent, we  can  put  here
$t=t_{0}$for a fixed $t_{0}\in J$, so  we  get  (2.5)  with  $C(\gamma
)={\bigl(}F^{-1}(t_{0};\gamma ){\bigr)}^{T}   G(\gamma
(t_{0}))F^{-1}(t_{0};\gamma )$, besides $C^{T}=C$ due to $G^{T}=$G. On the
contrary, if (2.5) is fulfilled, then (2.6) is  identically  valid  and,
consequently, (2.2) is also true. From here we conclude that  (2.2)  and
(2.6) are equivalent, i.e. propositions 2.1 and 2.2 are equivalent, the
former of which was already proved.\blacksquare

The arbitrariness in the choice of $F$ is described by the transformation
(1.5). As a consequence of this the matrix $C(\gamma )$ in (2.5) is also not
uniquely defined. Evidently, the  above  transformation leads to $C(\gamma
)\mapsto {\bigl(}D^{-1}(\gamma ){\bigr)}^{T}C(\gamma )D^{-1}(\gamma )$.
It is easy to check  the validity of the inverse to this statement, i.e. we
have
\[
F(s;\gamma )\mapsto D(\gamma )F(s;\gamma ) \Leftarrow \Leftrightarrow
C(\gamma )\mapsto {\bigl(}D^{-1}(\gamma ){\bigr)}^{T}C(\gamma
)D^{-1}(\gamma ).\qquad (2.7)
\]
 At this place naturally arises the question
when and under what conditions there exist bundle metrics (resp. transports
along paths) which are consistent with a given transport along paths (resp.
bundle metric). This question can be put in two variants: global, when
consistency along every path is investigated, and local, when consistency
along a given path is studied. The following propositions give a solution of
the above question from different viewpoints.

{\bf Proposition 2.3.} If an $L$-transport along paths is defined by the
matrix function $F$ through (1.4), then any consistent with it along $\gamma
:J\to B$ bundle metric along $\gamma $ has a matrix of the form
\[
 G(\gamma (s);\gamma )={\bigl(}F(s;\gamma ){\bigr)}^{T}C(\gamma
)F(s;\gamma ),\quad  s\in J,\ C^{T}=C,\qquad (2.8)
\]
 where $C(\gamma )$ is
nondegenerate, symmetric and depending on $\gamma $ matrix.

{\bf Proof.} This proposition is a corollary from proposition 2.2: in fact (2.8) is the general solution of (2.5) with respect to $G$ when $F$ is given.\blacksquare

Proposition 2.3 shows that locally, i.e. along a given path, every $L$-transport along paths defines (see $(2.8)) a$ unique class of consistent along that path with it metrics, uniquely defined only on the same path. The global variant of this problem is more difficult and will be treated below in proposition 2.6.

{\bf Proposition} ${\bf 2}{\bf .}{\bf 4}{\bf .} A$ necessary and sufficient condition for the
existence of globally (resp. locally) consistent with a given bundle metric $L$-transport along paths is the independence of the signature of the bundle metric, i.e. of the matrix $G(x)$, from the point of the base $B ($resp. the path of transport) at which it is evaluated, i.e. from $x\in B ($resp. $x\in \gamma (J))$.

 {\bf Remark.} If $p(x)$ and $q(x)$ are the numbers, respectively, of the
positive and negative eigenvalues of $G(x)$, the signature of $G(x)$ is
$s(x):=p(x)-q(x) ($see [6]). The proposition states that the bundle metric
admits a (globally or locally) consistent with it $L$-transport along paths
iff $s(x)=$const, which due to $p(x)+q(x)=n:=\dim(\pi ^{-1}(x)), x\in B$
is equivalent to $p(x)=$const and/or $q(x)=$const, i.e. to the independence
of the number of the positive (negative) eigenvalues of the metric from the
point at which they are evaluated.

{\bf Proof.} Let the bundle metric of $(E,\pi ,B)$ be (globally or locally)
consistent with an $L$-transport along paths described in some basis along
$\gamma :J\to B$ by the matrix $F_{0}(s;\gamma )$ through (1.4) (with
$F=F_{0})$. Then, by proposition 2.2, there exists a nondegenerate and
symmetric matrix $C_{0}(\gamma )$ such that ${\bigl(}F^{-1}_{0}(s;\gamma
){\bigr)}^{T}G(\gamma (s))F^{-1}_{0}(s;\gamma )= =C_{0}(\gamma ), s\in $J.
Because of $C^{T}_{0}(\gamma )=C_{0}(\gamma )$ there exists an orthogonal
matrix $D_{0}(\gamma )$ for which $D^{T}_{0}(\gamma )C_{0}(\gamma
)D_{0}(\gamma )$ is {\sl constant} diagonal matrix [6]. Then (see (1.5) and
(2.7)) the matrix $F(s;\gamma ):=D_{0}(\gamma )F_{0}(s;\gamma )$ describes
the same $L$-transport and due to the last equality it satisfies (2.5) with
$C(\gamma )=D^{T}_{0}(\gamma )C_{0}(\gamma )D_{0}(\gamma
)=$diag$(c_{1},\ldots  ,c_{n}), c_{1},\ldots  ,c_{n}\in  \in {\Bbb R}$. This
means that in (2.5) the matrix $F(s;\gamma )$ may be chosen in such a way as
$C(\gamma )$ to be a constant diagonal matrix. If such a choice is already
made, then, due to $(2.5), G(\gamma (s))$ with the help of a transformation
of a form $D^{T}G(\gamma (s))D$, where $D=F^{-1}(s;\gamma )$, can be
transformed into the diagonal matrix diag$(c_{1},\ldots  ,c_{n})$ which does
not depend either on $\gamma (s)$ or on $\gamma $. On the other hand, as
$G^{T}=G$, there exists an orthogonal matrix $D_{1}(\gamma (s))$ such that
\[
D^{T}_{1}(\gamma (s)) G(\gamma (s))D_{1}(\gamma (s))=diag(g_{1}(\gamma
(s)),\ldots  ,g_{n}(\gamma (s))),
\]
 where $g_{i}(\gamma (s))\neq 0,\infty , i=1,\ldots  ,n$ are the eigenvalues
of $G(\gamma (s)) [6]$. From the last two results, as a consequence of the
Jacobi-Sylvester inertia law [6], it follows that the number $p$ of the
positive and the number $q=n-p$ of the negative eigenvalues of $G(\gamma
(s))$, i.e. of the bundle metric, are equal, respectively, to the number of
positive and negative diagonal elements of the matrix $D^{T}G(\gamma
(s))D=$diag$(c_{1},\ldots  ,c_{n})=$const. Consequently the numbers $p$ and
$q(=n-p)$ do not depend on the point $\gamma (s), s\in J$ in the local case
(consistency along a path) or on the point of $B$ in the global case at which
they are calculated, i.e. the signature of the bundle metric is
$s=p-q=$const.

On the contrary, let in $(E,\pi ,B)$ be given a bundle metric $g$ whose
signature $s=p-q$, and consequently, the number $p$ of its positive and the
number $q=n-p$ of its negative eigenvalues, do not depend on the point $x$ at
which they are calculated $(x\in B$ in the global case or $x\in \gamma (J)$
in the local case for some path $\gamma :J\to B)$. Because of $G^{T}=G$ there
exists an orthogonal matrix $D_{1}(x)$ such that $D^{T}_{1}(x)G(x)
D_{1}(x)=$diag$(g_{1}(x),\ldots  ,g_{n}(x))$, where $g_{i}(x)\neq 0,\infty ,
i=1,\ldots  ,n$ are the eigenvalues of $G(x) [6]$. If we put
\[
D_{2}(x)
:=diag(\mid g_{1}(x)\mid ^{-1/2},\ldots  ,\mid g_{n}(x)\mid ^{-1/2})
\]
and $D(x):=D_{1}(x)D_{2}(x)$, then, as a result of the supposition made
above, we get $D^{T}(x)G(x)D(x)=$diag$(\epsilon _{1},\ldots ,\epsilon _{n})$,
where $p$ of the numbers $\epsilon _{1},\ldots ,\epsilon _{n}$are equal to
+1, the remaining $q=n-p$ ones being equal to -1.  Then from proposition 2.2
for $C(\gamma )=$diag$(\epsilon _{1},\ldots ,\epsilon _{n})=$const (see also
(2.5)) and (1.4) for $F(s;\gamma )={\bigl(}D(\gamma (s)){\bigr)}^{-1},
s\in J$ for every path $\gamma :J\to B ($global case) or for some path
$\gamma :J\to B ($local case), we conclude that the given bundle metric is,
respectively, globally or locally consistent with the $L$-transport along
paths defined by the matrix (1.4) with the above definition of F.\blacksquare

The next proposition describes the general form and the way of construction of $L$-transports along paths consistent with a given bundle metric admitting such transports along paths.

{\bf Proposition 2.5.} Let in $(E,\pi ,B)$ be given a bundle metric whose
signature does not depend on the point $\gamma (s)$ at which it is evaluated
for every (resp. some) path $\gamma :J\to $B. Let there be chosen bases
$\{e_{i}(s)\}, s\in J$ along $\gamma $ in such a way that the first $p$
eigenvalues of the matrix $G(\gamma (s))$, defining the metric in them (see
(2.1)), be positive. Then one $L$-transport along paths is consistent with
this bundle metric along every (resp. a given) path $\gamma $ if and only if
some of the defining its matrix (1.4) matrix functions $F$ has the form
\[
 F(s;\gamma )=Y(\gamma )Z(s;\gamma ){\bigl(}D(\gamma (s))
{\bigr)}^{-1},\quad  s\in J\qquad (2.9)
\]
 for every (resp. the given) path $\gamma $. In
this equality:  $Y(\gamma )$  is  $n n$ nondegenerate matrix function of
$\gamma ; Z(s;\gamma )$  is  a  pseudo-orthogonal matrix of type $(p,q),
q=n-p$, i.e. $Z(s;\gamma )\in O(p,q)$, or
\[
{\bigl(} Z(s;\gamma ){\bigr)}^{T} G_{p,q} Z(s;\gamma )
=G_{p,q}
:=diag(\underbrace{1,\ldots,1}_{p\ times},
       \underbrace{-1,\ldots,-1}_{q\ times})
\qquad (2.10)
\]
 and $D(\gamma (s))$ is a fixed (orthogonal) matrix such that
\[
{\bigl(}D(\gamma (s)){\bigr)}^{T}G(\gamma (s))D(\gamma
(s))=G_{p,q}.\qquad (2.11)
\]

 {\bf Remark.} The case when in some basis not all
of the first $p$ eigenvalues of the bundle metric are positive is obtained
from the above one by transformation (renumbering) of the basis
$\{e_{i}(s)\}$.

{\bf Proof.} To prove the necessity we have to solve the equation (2.5),
when $G$ and $C$ are given, with respect to $F(s;\gamma )$. From $G^{T}=G$,
the choice of $\{e_{i}(s)\}$ and the independence of $p$ and $q=n-p$ from
$\gamma (s) ($because of $p-q=s=$const) follows the existence of a satisfying
(2.11) matrix $D(\gamma (s))$: e.g. we can put $D(\gamma (s))=D_{1}(\gamma
(s))D_{2}(\gamma (s))$, where $D_{1}$and $D_{2}$were defined in the proof of
proposition 2.4.

Let $F(s;\gamma )=:F_{1}(s;\gamma ){\bigl(}D(\gamma (s)){\bigr)}^{-1}$.
From equations (2.5) and (2.11), we get ${\bigl(}F^{-1}_{1}(s;\gamma
){\bigr)}^{T}G_{p,q}F^{-1}_{1}(s;\gamma )=C(\gamma )$. Putting here
$F_{1}(s;\gamma )=:Y(\gamma )Z(s;\gamma )$, where $Y(\gamma )$ is arbitrary
nondegenerate matrix for which ${\bigl(}Y(\gamma ){\bigr)}^{T}C(\gamma )
Y(\gamma )=G_{p,q}$,  we see that $Z(s;\gamma )$ satisfies (2.10). (The
existence of $Y(\gamma )$ follows from (2.5): from it and the law of
Jacobi-Sylvester it follows that $C(\gamma )$ has $p$ positive eigenvalues,
due to which the needed matrix $Y(\gamma )$, which is orthogonal, exists
[6].) All that proves that $F(s;\gamma )$ has the form (2.9) if the
considered $L$-transport along paths is consistent along every (resp. some)
path $\gamma $ with the given bundle metric.

On the contrary, the sufficiency of the proposition is almost evident: if (2.9) is valid, then it is easy to check the validity of (2.5) for $C(\gamma )={\bigl(}Y^{-1}(\gamma ){\bigr)}^{T}G_{p,q}Y^{-1}(\gamma )$ and according to proposition 2.2 the $L$-transport along paths and the bundle metric are consistent along every (resp. some) path $\gamma .\blacksquare $

{\bf Corollary 2.1.} For a given $L$-transport along paths $L$ there locally exists a consistent with it along $\gamma $ bundle metric if and only if along $\gamma $ exists a basis in which the defining it by (1.4) matrix $F(s;\gamma )$ has the form (2.9), in which $Y(\gamma )$ and $D(s;\gamma )$ are arbitrary $n n$ nondegenerate matrices and $Z^{T}(s;\gamma )G_{p,q}Z(s;\gamma )=G_{p,q}$for some $p,q\ge 0, p+q=n$ which may depend on $\gamma $.

{\bf Proof.} If along $\gamma $ there exists a consistent with $L$ bundle
metric $g$, then the expansion (2.9) follows from proposition 2.5.
Conversely, if (2.9) is valid, then substituting (2.9) into (2.8), we get, in
accordance with proposition $2.3, a$ class of consistent with $L$ along
$\gamma $ metrics defined along $\gamma $. In particular, choosing $Y$ and
$C$ in such a way that $Y^{T}CY=G_{p,q}$,  we obtain
\[
G(\gamma (s);\gamma )=((D(\gamma (s)))^{-1})G_{p,q}(D(\gamma (s)))^{-1},
\]
 where $p$ and $q$ may depend on $\gamma .\blacksquare $

Now we shall go back to the question of the global existence of a bundle metric globally consistent with a given $L$-transport along paths (see the comment after the proof of proposition 2.3).

{\bf Proposition} ${\bf 2}{\bf .}{\bf 6}{\bf .} A$ necessary and sufficient
condition for the existence of globally consistent with a given $L$-transport
along paths bundle metrics is the existence of a local basis along any path
$\gamma :J\to B$ in which the matrix $F(s;\gamma )$, defining this transport
along $\gamma :J\to B$ by (1.4), has the form (2.9) in which: $Y(\gamma )$
and $D(\gamma (s))$ are arbitrary nondegenerate $n n$ matrices and
$Z(s;\gamma )$ is a pseudo-orthogonal matrix of type $(p,q)$, with  $p+q=$n.
Besides, if for a given $L$-transport along paths are fulfilled these
conditions, then in the above mentioned bases all globally consistent with
the $L$-transport along paths bundle metrics are defined by the matrix
\[
 G(\gamma (s))={\bigl(}(D(\gamma (s)))^{-1}{\bigr)}^{T}G_{p,q}
{\bigl(}D(\gamma (s)){\bigr)}^{-1},\qquad (2.12)
\]
 which depends only on the point
$\gamma (s)$, but not on the path $\gamma $.

{\bf Proof.} Let there be given an $L$-transport along paths for which there
exists a globally consistent with it bundle metric. Now we shall prove that
under this condition is fulfilled an equality like (2.9). In arbitrary bases
along every path $\gamma $ this $L$-transport along $\gamma $ defines a class
of consistent with it bundle metrics along $\gamma $ which are defined by
(2.8). Let in (2.8) the matrix $C(\gamma )$ have $k(\gamma )$ positive and
$l(\gamma )=n-k(\gamma )$ negative eigenvalues. Then (see [6]) there exists a
matrix $Y(\gamma )$ such that $Y^{T}(\gamma )C(\gamma )Y(\gamma )=G_{k(\gamma
),l(\gamma )}$, as a consequence of which (2.8) can be represented as
\[
G(\gamma (s);\gamma )=F^{T}(s;\gamma ){\bigl(}Y^{-1}(\gamma )
{\bigr)}^{T}G_{k(\gamma ),l(\gamma )}Y^{-1}(\gamma )F(s;\gamma ).
\]
If the consistent with an $L$-transport along paths metric, which by
assumption exists, is described in some basis with the matrix $G_{0}(x)$
having for every $x\in B p$ positive and $q=n-p$ negative eigenvalues, then
by proposition 2.3 along $\gamma $ this bundle metric belongs to the above
class of bundle metrics and, hence, there exists a matrix $Y_{0}(\gamma )
($or the corresponding matrix $C_{0}(\gamma )$; see (2.8)) such that
\[
G_{0}(\gamma (s);\gamma )=F^{T}(s;\gamma ){\bigl(}Y^{-1}_{0}(\gamma
){\bigr)}^{T}G_{k_{0}}Y^{-1}_{0}(\gamma )F(s;\gamma ).
\]

 From here, due to the Jacobi-Sylvester law [6], it follows that
$k_{0}(\gamma )=p=$const and $l_{0}(\gamma )=n-p=$const. If we put
$F(s;\gamma )=:Y_{0}(\gamma )   Z_{0}(s;\gamma )D^{-1}_{0}(\gamma (s))$,
where $D_{0}(x)$ is a matrix for which $D^{T}_{0}(x)G_{0}(x)
D_{0}(x):=G_{p,q}, x\in B [6]$, then from the same equality, we see that
\[
G_{0}(\gamma (s))={\bigl(}(D_{0}(\gamma (s)))^{-1}
{\bigr)}^{T}Z^{T}_{0}(s;\gamma )G_{p,q}Z_{0}(s;\gamma ){\bigl(}D_{0}(\gamma
(s)){\bigr)}^{-1}.
\]
 From here, as a consequence of the definition of $D_{0}$, we conclude that
$Z_{0}(s;\gamma )$ is a pseudo-orthogonal matrix of type $(p,q) ($see
(2.10)). This result proves the existence of a representation like (2.9) for
the considered $L$-transport along paths.

The sufficiency of the proposition is almost evident. In fact, if (2.9) is
valid for one $L$-transport along paths, then substituting (2.9) into (2.8),
we see the class of consistent with it along $\gamma $ bundle metrics to be
defined by
\[
G(\gamma (s);\gamma )={\bigl(}(D(\gamma (s)))^{-1}{\bigr)}^{T}
Z^{T} (s;\gamma )Y^{T}(\gamma )C(\gamma )Y(\gamma )Z(s;\gamma )
{\bigl(}D(\gamma (s)){\bigr)}^{-1},
\]
 where $C(\gamma )$ is nondegenerate, symmetric $(C^{T}=C)$ matrix of type $n
n$, which plays a role of a parameter whose change describes the considered
class of bundle metrics along $\gamma $. If $p$ and $q$ are arbitrary
nonnegative integers and $p+q=n$, then, due to the last equality and (2.10),
the choice
\[
C(\gamma )={\bigl(}Y^{-1}(\gamma ){\bigr)}^{T}G_{p,q}Y^{-1}(\gamma
)=C^{T}(\gamma )
\]
 defines a set of $n+1 ($resp. for $p=0, \ldots  . ,p=n; q=n-p)$ bundle
metrics along $\gamma $, given by the equality
\[
 G(\gamma (s);\gamma ;p,q)={\bigl(}D^{-1}(\gamma (s))
{\bigr)}^{T}G_{p,q}D^{-1}(\gamma (s))\qquad (2.12^\prime )
\]
and, evidently, depending only on the point $\gamma (s)$, but not on the map
$\gamma :J\to B$. Consequently these bundle metrics are globally consistent
with the given $L$-transport along paths.\blacksquare

\medskip
\medskip
 {\bf 3. THE CASE OF  GENERATED BY DERIVATIONS OF\\
	 TENSOR	ALGEBRAS LINEAR TRANSPORTS OF VECTORS}

\medskip
In this section we shall concentrate our attention on the tangent bundle $(T(M),\pi ,M)$ of a given manifold M. In particular we are going to consider in it some special questions concerning the consistency of bundle metrics and $S$-transports (linear transports along paths generated by derivations of tensor algebras [8]) of vectors along (smooth) paths. In this case the bundle metric $g$ is a nondegenerate symmetric section (tensor field) of type (0,2). The bundle metrics in this fibre bundle are used to be called simply metrics [3,7], because of which the adjective "bundle", as applied to metric(s), will be omitted till the end of the present section.

The details concerning $S$-transports along paths and their properties can be found in [8].

In accordance with definition 1.1 the metric $g$ and the $S$-transport
$S^{\gamma }$along $\gamma :J\to M$ are (locally) consistent along $\gamma $
if
\[
g_{\gamma (s)}(A_{0},B_{0})=g_{\gamma (t)}(S^{\gamma }_{s\to
t}A_{0},S^{\gamma }_{s\to t}B_{0}), s,t\in J\qquad (3.1)
\]
 for every
$A_{0},B_{0}\in T_{\gamma (s)}(M)$. The metric and the $S$-transport are
(globally) consistent if this equality is fulfilled for every path $\gamma $.

 {\bf Proposition 3.1.} The $C^{1}$metric $g$ and the $S$-transport $S$ are
consistent (resp. along a given path $\gamma )$ iff the equality
\[
 g_{\gamma (t)}=S^{\gamma }_{s\to t}g_{\gamma (s)},\quad s,t\in J\qquad (3.2)
\]
 is valid for every (resp. the given) path $\gamma $.

{\bf Proof.} Using the contraction operator $C^{1}_{1}$on the first
superscript and the first subscript, we find
\[
 g(A,B)=(C^{1}_{1})^{2}(A\otimes g\otimes B), A,B\in Sec(T(M),\pi
,M).\qquad (3.3)
\]
If we apply $S^{\gamma }_{t\to s}$to (3.1) and take into
account (3.3), the properties of the $S$-transports [8], and the
arbitrariness of $A_{0}$and $B_{0}$, we get that (3.1) is equivalent to
$g_{\gamma (s)}=S^{\gamma }_{t\to s}g_{\gamma (t)}$which, due to the
arbitrariness of $s,t\in J$, is equivalent to (3.2).\blacksquare

Let us note that by definition 5.1 from [8] the equality (3.2) means that the metric $g$ is $(S-)$transported along $\gamma $ section, hence proposition 3.1 is equivalent to

{\bf Proposition 3.1'.} The $C^{1}$metric $g$ and an $S$-transport $S$ are consistent (resp. along a given path $\gamma )$ iff $g$ is $S$-transported
along every (resp. the given) path $\gamma $ section.

{\bf Proposition 3.2.} The $S$-transport $S^{\gamma }$along $\gamma :J\to M$
is consistent along $\gamma $ with the $C^{1}$metric $g$ iff
\[
({\cal D}^{\gamma }g)(\gamma (s))={\cal D}^{\gamma }_{s}g=0, s\in
J,\qquad (3.4)
\]
 where ${\cal D}^{\gamma }$is the defined by or defining
$S^{\gamma }$derivation of the tensor algebra over $\gamma (J)$ according to
[8], proposition 8.

{\bf Proof.} This proposition is a consequence of proposition $3.1^\prime $ and proposition 5.3. of [9].\blacksquare

 The equality (3.4) is useful and effective tool for practical check of the
question of consistency of $S$-transports and $C^{1}$metrics. A typical
example for this is the Riemannian parallel transport (in (pseudo-)Riemannian
manifolds) which is defined from the Levi-Cevita connection and for which
(3.4) is identically satisfied $[3, 4]$: in this case $\nabla
^{\{\}}_{X}g=0$, where $\nabla ^{\{\}}_{X}$is the defined by the Christoffel
symbols from $g$ covariant differentiation along the vector field $X$ and as
now ${\cal D}^{\gamma }=\nabla _{\cdot }$,   being the tangent to $\gamma $
vector field $(cf. [2]$, Sect. 5), the $eq. (3.4)$ holds identically. In
(pseudo-) Riemannian manifolds the connection for which $\nabla _{X}g=0$ is
called metric preserving as it preserves the scalar product of the (parallel)
transported with its help vectors along any path [10].

For some purposes it is convenient to write the conditions (3.2) and (3.4) in local coordinates. For this goal along a path $\gamma :J\to M$ are introduced bases $\{E_{i}\mid _{\gamma (s)}\}$ in the tangent spaces $T_{\gamma (s)}(M)$ to $M$ at $\gamma (s), s\in J$ in which the metric along $\gamma $ is given from the matrix $G(\gamma (s))$ with components (2.1), and the derivation ${\cal D}^{\gamma }$and the $S$-transport $S^{\gamma }$along $\gamma $ are described with the help of the matrix $\Gamma _{\gamma }(s):= \Gamma ^{i}_{.j}(s;\gamma ) $ of the coefficients of ${\cal D}^{\gamma }$and $S$ in $\{E_{i}\} ($see [2,8]).

{\bf Proposition 3.3.} The $C^{1}$metric $g$ and the $S$-transport $S$ are
consistent (resp. along a given path $\gamma )$ iff one of the following two
equalities is valid
\[
 G(\gamma (t))=Y(t,s;\Gamma _{\gamma })G(\gamma (s))Y^{T}(t,s;\Gamma
_{\gamma }), s,t\in J,\qquad (3.5)
\]
\[
\frac{dG(\gamma(s)}{ds} - G(\gamma (s))\Gamma _{\gamma }(s) - \Gamma _{\gamma
}(s)G(\gamma (s))=0,\quad s\in J,\qquad (3.6)
\]
 where $Y$ is the solution of
the initial-value problem (3.9) of [8], for every (resp. the given) path
$\gamma $. The equalities (3.5) and (3.6) are equivalent between each other
as well as, respectively, to (3.2) and (3.4). Besides, (3.5) is the general
solution of (3.6) with respect to G.

{\bf Proof.} With the help of $(3.4), (2.16)$ and $[8], eq. (3.8)$ it is easy to check that in the chosen bases (3.5) and (3.6) are component form in matrix notation of (3.2) and (3.4) respectively and, consequently, (3.5) and (3.6) are equivalent, respectively, to (3.2) and (3.4). From here and propositions 3.1 and 3.2 follows the first part of the proposition. From the same propositions follows also the equivalence of (3.2) and (3.4), and hence the equivalence of (3.5) and (3.6). If we look on (3.6) as an equation with respect to $G ($resp. along $\gamma )$, then from the definition of $Y$ it is clear that its general solution is given by (3.5) in which $s$ has to be fixed and $G(\gamma (s))$ must be replaced with an arbitrary constant matrix $C ($resp. a matrix function $C(\gamma )$ of $\gamma ).\blacksquare $

{\bf Proposition 3.4.} If in a given basis one $S$-transport is defined by
the matrix $\Gamma _{\gamma }$of its coefficients and the metric $g$ - by the
matrix $G$, then the $S$-transport and the metric are consistent (resp. along
a given path $\gamma )$ iff there exists a nondegenerate symmetric matrix
$C(\gamma )$ such that
\[
 Y(s_{0},s;\Gamma _{\gamma })G(\gamma (s))Y(s,s_{0};-\Gamma _{\gamma
})=C(\gamma )=C^{T}(\gamma )\qquad (3.7)
\]
 for $Y$ defined in $[8], eq. (3.9),
s\in J$, fixed $s_{0}\in J$ and every (resp. the given) path $\gamma $.

{\bf Proof.} According to the proof of proposition 3.4 of [8] the
$S$-transport is uniquely defined by the matrices
\[
F(s;\gamma )=Y(s_{0},s;-\Gamma _{\gamma })=[Y(s,s_{0};-\Gamma _{\gamma
})]^{-1}, s\in J\qquad (3.8)
\]
 through (1.4). Substituting (3.8) into (2.5)
(with the help of the properties of $Y)$, we get (3.7), which shows that
proposition 3.4 is a special case of proposition 2.2.\blacksquare

{ \bf Proposition 3.5.} If an $S$-transport is fixed in a given basis through
the matrix $\Gamma _{\gamma }$, then any consistent along $\gamma :J\to M$
with it metric is obtained in the same basis by the formula
\[
 G(\gamma (s);\gamma )=Y(s,s_{0};\Gamma _{\gamma })C(\gamma
)Y(s_{0},s;-\Gamma _{\gamma }),\qquad (3.9)
\]
 where $Y$ is defined in $[8],
eq. (3.9), s_{0}\in J$ is fixed and $C(\gamma )$ is nondegenerate symmetric
matrix. A necessary condition for the obtained by this formula metrics to be
globally consistent with the initial $S$-transport is $C(\gamma )$ to be
independent of $\gamma $, i.e.  $C(\gamma )=$const.

{\bf Proof.} The first part of the proposition follows, analogously to the
proof of proposition 3.4., from proposition 2.3 for $F(s;\gamma )=
=Y(s_{0},s;-\Gamma _{\gamma })$ and the fact that (3.9) is the general
solution of (3.6) with respect to $G ($see proposition 3.3). The second part
of the proposition follows from the circumstance that if in (3.9) we put
$s=s_{0}$, we get
\[
  G(\gamma (s_{0});\gamma )=C(\gamma ).\qquad (3.10)
\]
 Hence, if the given by
(3.9) metrics and the defined by $\Gamma _{\gamma }S$-transport are
(globally) consistent, then $C(\gamma )$ will not depend on $\gamma $ as
$G(\gamma (s_{0});\gamma )$ is simply the value of the matrix $G(x)$, which
represents $g$ in the given basis, at the point $\gamma (s_{0})$, but,
evidently, $G(\gamma (s_{0}))$ does not depend on $\gamma .\blacksquare $

{ \bf Proposition 3.6.} Let $\gamma :J\to M$ be $a C^{1}$path and $g$ be $a
C^{1}$metric (resp. along $\gamma )$ which in a basis $\{E_{i}\mid _{\gamma
(s)}\}$ along $\gamma $ is represented by the matrix $G(x), x\in \gamma (J)$,
whose signature does not depend on the point x. Let $D(x)$ be an orthogonal
matrix such that [6]
\[
 D^{T}(x)G(x)D(x)=\tilde{G}(x):=diag(g_{1}(x),\ldots  ,g_{n}(x)),
D^{T}=D^{-1},\qquad (3.11)
\]
 where $g_{i}(x)\neq 0,\infty , i=1,\ldots  ,n$
are the eigenvalues of $G(x)$ and
\[
 K(\gamma (s))
:=D^{T}(\gamma (s))\cdot \frac{dG(\gamma(s)}{ds} \cdot
D(\gamma (s))=K^{T}(\gamma (s)),\quad s\in J.   \qquad (3.12)
\]
 Then in a considered
basis the matrix $\Gamma _{\gamma }(s)$ of the coefficients of all consistent
(resp.  along $\gamma )$ with the metric $g S$-transports has the form
\[
 \Gamma _{\gamma }(s)=D(x)[P(x)+Q(x)\tilde{G}(x)+R(x)]D^{T}(x)
\]
\[
 =D(x)[P(x)+R(x)]D^{T}(x)+D(x)Q(x)D^{T}(x)G(x), x=\gamma (s), s\in \gamma
(J),  \qquad (3.13)
\]
 where the matrices $P, Q, R$ have the following symmetries:
\[
 P^{T}=P, Q^{T}=-Q, R^{T}=-R\qquad (3.14a)
\]
 and their components are:
\[
 Q_{ij}(x) =0 \quad
R_{ij}(x) = K_{ij}(x)/2g_i(x)
\qquad  for\ g_{i}(x)+g_{j}(x)=0,\qquad (3.14b)
\]
\[
P_{ij} = K_{ij}(x)/(g_i(x)+g_j(x)) \quad
R_{ij}(x) = 0\qquad  for\ g_{i}(x)+g_{j}(x)\neq 0,\qquad (3.14c)
\]
 the remaining of which can be chosen arbitrarily (only if (3.14a) are
satisfied).

{\bf Proof.} In fact we have to prove that (3.13) is the general solution of
(3.6) along $\gamma $ with respect to $\Gamma _{\gamma }$when $G$ is given.

Multiplying (3.6) on the left by $D^{T}(\gamma (s))$ and on the right by
$D(\gamma (s))$, and using (3.11), we get
\[
 \tilde{G}\tilde{\Gamma }+\tilde{\Gamma }\tilde{G}=K,\qquad (3.15)
\]
 where  $\tilde{\Gamma }:=D^{T}\Gamma _{\gamma }D$ and here, as well as below
in this proof, for brevity we omit everywhere the arguments $s$ and $\gamma
(s), s\in $J. When written in component form (see (3.11)) this equation will
be equivalent to
\[
  g_{i}\tilde{\Gamma }_{ij}+g_{j}\tilde{\Gamma }_{ji}=K_{ij}.\qquad
(3.15^\prime )
\]
 (Do not sum here over $i$ and $j$!)

Let us consider the pairs $(i,j)$ for which $g_{i}+g_{j}=0$. Then $(3.15^\prime )$ reduces to $g_{i}(\tilde{\Gamma }_{ij}-\tilde{\Gamma }_{ji})=K_{ij}$, hence we can find only the antisymmetric part of the element $\tilde{\Gamma }_{ij}$. So, using the identity $\tilde{\Gamma }_{ij}=$ $(\Gamma _{ij}+\tilde{\Gamma }_{ji})+$ $(\Gamma _{ij}-\tilde{\Gamma }_{ji})$, we get $\tilde{\Gamma }_{ij}=P_{ij}+R_{ij}+Q_{ij}g_{j}$, where $R_{ij}:= :=K_{ij}/2g_{i}=:-R_{ji}, Q_{ij}:=-Q_{ji}:=0$ and the quantities $P_{ij}:=P_{ij}$are arbitrary.

Now we shall consider the pairs $(i,j)$ for which $g_{i}+g_{j}\neq 0$. In this case we define the quantities $P_{ij}=P_{ji}$, i.e. the remaining components of the matrix $P=P^{T}$, as the symmetric solution of the equation (3.15 ), i.e. $g_{i}P_{ij}+g_{j}P_{ji}=K_{ij}$, and hence $P_{ij}=K_{ij}/(g_{i}+g_{j})= =P_{ji}$. Then, if we put $\tilde{\Gamma }_{ij}=:P_{ij}+Q_{ij}g_{j}+R_{ij}, R_{ij}:=-R_{ji}:=0$, we see that (3.15 ) reduces to $Q_{ij}+Q_{ji}=0$, i.e. the only restriction on the quantities $Q_{ij}$is their antisymmetry.

Thus we proved that $\tilde{\Gamma }=P+Q\tilde{G}+R$, where $P, Q$ and $R$ are defined by (3.14), is the general solution of (3.15) with respect to $\tilde{\Gamma }:=D^{T}\Gamma _{\gamma }$D. From here it follows that $\Gamma _{\gamma }=(D^{T})^{-1}\tilde{\Gamma }D^{-1}=D\tilde{\Gamma }D^{T}$is given by (3.13) and
it is the general solution of (3.6) with respect to $\Gamma _{\gamma }.\blacksquare $

\medskip
\medskip
 {\bf 4. CONSISTENCY BETWEEN LINEAR TRANSPORTS ALONG\\
	PATHS AND HERMITIAN BUNDLE METRICS}

\medskip
In this section with $(E,\pi ,B)$ is denoted an arbitrary complex vector bundle [1,11].

By a Hermitian bundle metric $g$ in $(E,\pi ,B)$ we understand (see [1] and
$[11], ch. I, \S8) a$ map $g:x\mapsto g_{x}, x\in B$, where the maps
\[
  g_{x}:\pi ^{-1}(x) \pi ^{-1}(x)\to {\Bbb C}\qquad (4.1)
\]
 have the properties:  $1\frac{1}{2}$  linearity, nondegeneracy and
Hermiticity, i.e.
\[
 g_{x}(u+v,w)=g_{x}(u,w)+g_{x}(v,w),\quad u,v,w\in \pi ^{-1}(x),\qquad (4.2a)
\]
\[
g_{x}(\lambda u,\mu v)=\lambda g_{x}(u,v),\quad \lambda ,\mu \in {\Bbb C},
\ u,v\in \pi ^{-1}(x),\qquad (4.2b)
\]
\[
 g_{x}(u,v)\neq 0\quad for\ u,v\neq {\bf O}\in \pi ^{-1}(x),\qquad (4.2c)
\]
\[
g_{x}(u,v)=\overline{g_{x}(v,u)} ,  \qquad (4.2d)
\]
 where the bar over a complex number or a matrix means complex conjugation.
In this definition we neglect the usual [11] condition $g_{x}(u,u)>0$ for
$u\neq {\bf 0}$ as insignificant for the following investigation.

The studied below $L$-transports along paths in the complex vector bundle $(E,\pi ,B)$ are supposed to be ${\Bbb C}$-linear $(cf. [2]$, Sect. 2 or [14], Subsect. 2.3).

Evidently, in the real case, i.e. when ${\Bbb C}$ is replaced with ${\Bbb R}, a$ Hermitian bundle metric reduces to the defined in section 1 real bundle metric. Therefore any result for Hermitian bundle metrics is also valid for real bundle metrics. The contrary is mutatis mutandis true, i.e. with certain changes in the formulations of the definitions and propositions of sections 1 and 2 they remain valid also in the Hermitian case. The present section is intended for the description of these changes.

 The basic definition 1.1 now reads.

{\bf Definition 4.1.} The linear transport along paths $L$ and the Hermitian
bundle metric $g$ are called consistent (resp. along the path $\gamma )$ if
$L$ preserve the Hermitian scalar products of the vectors along every (resp.
along the given) path $\gamma :J\to B$, i.e. if
\[
 g_{\gamma (s)}=g_{\gamma (t)}\circ (L^{\gamma }_{s\to t} L^{\gamma
}_{s\to t}),\quad s,t\in J\qquad (4.3)
\]
for every (resp. for the given) path $\gamma $.

Let $\gamma :J\to B$ and in every fibre $\pi ^{-1}(\gamma (s)), s\in J$ be
fixed a (complex) basis $\{e_{i}(s): i=1,\ldots  ,n=\dim_{{\Bbb C}}(\pi
^{-1}(x)), x\in B\}$. In them the $L$-transport $L^{\gamma }_{s\to t}$along
$\gamma $ from $s$ to $t$ is uniquely defined by its, generally complex,
matrix $H(t,s;\gamma )= H^{\hbox{i.}}_{.j}(t,s;\gamma )  ($see [2]); the
metric $g$ at the point $\gamma (s)$ is defined by the matrix $G(\gamma
(s)):= (g_{\gamma (s)})_{ij} , i,j=1,\ldots  ,n=\dim_{{\Bbb C}}(\pi
^{-1}(\gamma (s))$ with the defined by (2.1) elements. In terms of $G(\gamma
(s))$ the nondegeneracy and Hermiticity of $g$ mean, respectively,
$\det(G(\gamma (s)))\neq 0,\infty $ and $G^{*}(\gamma (s))=G(\gamma (s))$,
where * means Hermitian conjugation of matrices $(G^{*}:=$ $^{T}=($ $)^{T}$;
see [6]), i.e. $G$ is a nondegenerate Hermitian matrix function. This follows
from the fact that if $u=u^{i}e_{i}(s)$ and $v=v^{i}e_{i}(s)$, then due to
(4.2), we have
\[
u^{i}\overline{v^j} g_{\gamma (s)}(e_{i}(s),e_{j}(s))
=g_{\gamma (s)}(u,v)
= \overline{g_{\gamma (s)}(v,u) } =
\]
\[
 = \overline {
v^j \overline{u^i} g_{\gamma(s)} (e_{j}(s))
}
,\qquad (4.4)
\]
 i.e.  $\det[(g_{\gamma (s)})_{ij}] \neq 0$ and
$(g_{\gamma (s)})_{ij} = \overline{ ( g_{\gamma(s)} )_{ji} }$.

Below we present the analogs of propositions $2.1-2.6$, respectively, with numbers $4.1-4.6$ for Hermitian bundle metrics.

{\bf Proposition} ${\bf 4}{\bf .}{\bf 1}{\bf .} A$ necessary and sufficient
condition for global (resp. local) consistency of a Hermitian bundle metric
$g$ and an $L$-transport along paths $L$ is the equality
\[
G^{T}(\gamma (s))=H^{*}(t,s;\gamma )G^{T}(\gamma (t))H(t,s;\gamma ),
s,t\in J\qquad (4.5)
\]
 for every (resp. for a given) path $\gamma :J\to $B.

{\bf Proof.} Applying (4.3) to $(u,v), u,v\in B$, using (4.4) and the
arbitrariness of $(u,v)$, we get
\[
 (g_{\gamma (s)})_{ij}
=(g_{\gamma (t)})_{kl}H^{\hbox{k.}}_{.i}(t,s;\gamma )
\overline{ H^{\hbox{i.}}_{.j}(t,s;\gamma ) }
\qquad (4.5^\prime )
\]
from where taking a complex conjugate and taking into account $G^{T}=$ , due
to $G^{*}=G$, we obtain the component form of (4.5).\blacksquare

{\bf Proposition 4.2.} If the nondegenerate (complex) matrix function
$F(s;\gamma )$ defines an $L$-transport along paths $L$ by (1.4) then a
necessary and sufficient condition for the global (resp. local) consistency
of $L$ with a Hermitian bundle metric $g$ is the existence of a nondegenerate
Hermitian matrix function $C$ of $\gamma $, such that
\[
 {\bigl(}F^{-1}(s;\gamma ){\bigr)}^{*}G^{T}(\gamma (s))F^{-1}(s;\gamma
)=C(\gamma ),\quad  s\in J,\ C^{*}=C\qquad (4.6)
\]
for every (resp. a given) path $\gamma $.

{\bf Proof.} The proof of this proposition is an exact copy  of
the proof of proposition 2.2 and it  reduces  to  the  substitution  of (1.4)
into (4.5) and the separation in the  obtained  equality  the terms at
$\gamma (s)$ to the left of the equality sign and those at $\gamma (t)$  - to
its right side.\blacksquare

 As is known [2] the function $F$, appearing in (4.6), is defined up to the
transformation (1.5), so the function $C$ in (4.6) is not uniquely defined
and analogously to (2.7) one finds the implication
\[
 F(s;\gamma )\mapsto D(\gamma )F(s;\gamma ) \Leftrightarrow
C(\gamma )\mapsto {\bigl(}D^{-1}(\gamma ){\bigr)}^{*}C(\gamma
)D^{-1}(\gamma ).\qquad (4.7)
\]

 {\bf Proposition 4.3.} If an $L$-transport
along paths is defined by the matrix function $F$ through (1.4), then any
consistent with it along $\gamma :J\to B$ bundle Hermitian metric along
$\gamma $ is given by
\[
G^{T}(\gamma (s);\gamma )={\bigl(}F(s;\gamma ){\bigr)}^{*}C(\gamma
)F(s;\gamma ), s\in J, C^{*}=C,\qquad (4.8)
\]
 where $C$ is nondegenerate
Hermitian matrix function of $\gamma $.

{\bf Proof.} Solving (4.6) with respect to $G(\gamma (s))$, when $F$ is fixed and $C$ is arbitrary, we get (4.8).\blacksquare

{\bf Proposition} ${\bf 4}{\bf .}{\bf 4}{\bf .} A$ necessary and sufficient condition for the existence of globally (resp. locally) consistent with a given Hermitian bundle metric $L$-transport along paths is the independence of the signature (and consequently, of the number of positive (and/or negative) eigenvalues) of that metric, i.e. of the matrix $G(x)$, from the point of the manifold $M ($resp. of the path) at which it is evaluated, i.e. of $x\in B ($resp. $x\in \gamma (J))$.

{\bf Proof.} The proof of this proposition in form coincides with the one of
proposition 2.4 and as the latter is simple but long enough here we shall
present only the changes which must be done in it for obtaining the needed
proof in the Hermitian case.

1. The matrix $C_{0}(\gamma )$ is nondegenerate, Hermitian and such that
${\bigl(}F^{-1}_{0}(s;\gamma ){\bigr)}$
 $\times^{*}G^{T}(\gamma (s))F^{-1}_{0}(s;\gamma )
=C_{0}(\gamma )=C^{*}_{0}(\gamma )$.

2. The matrix $D_{0}(\gamma )$ is unitary $(D^{-1}_{0}=D^{*})$ and
satisfies the equation
$C(\gamma )=$$D^{*}_{0}(\gamma )C_{0}(\gamma )D_{0}(\gamma
)=diag(c_{1},\ldots  ,c_{n}), c_{1},\ldots  ,c_{n}\in {\Bbb C}$,
$n=\dim_{{\Bbb C}}(\pi ^{-1}(x))$, $x\in B$.

3. Further in the proof the matrix functions $D$ and $D_{1}$are unitary (not orthogonal$), D^{T}$and $D^{T}_{1}$must be replaced, respectively, with $D^{*}$and $D^{*}_{1}($due to $G^{*}=G$, but not $G^{T}=G)$ and one must have in mind that the eigenvalues of the Hermitian matrices (in this case $C_{0}$and $G)$ are real (e.g. $c_{1},\ldots  ,c_{n}\in {\Bbb R}\subset {\Bbb C}) [6].\blacksquare $

{\bf Proposition 4.5.} Let in $(E,\pi ,B)$ be given a Hermitian bundle
metric $g$ the signature of which is independent of the point $\gamma (s)$ at
which it is evaluated for every (resp. some) path $\gamma :J\to $B. Let there
be chosen bases $\{e_{i}(s)\}, s\in J$ along $\gamma $ such that the first
$p$ eigenvalues of the Hermitian matrix $G(\gamma (s))$, defining the metric
in them (see (2.1)), be positive. Then one $L$-transport along paths is
consistent with this Hermitian metric if and only if some of the defining it
by (1.4) matrix functions $F$ has the form
\[
 F(s;\gamma )=Y(\gamma )Z(s;\gamma ){\bigl(}D(\gamma (s))
{\bigr)}^{-1},\quad  s\in J\qquad (4.9)
\]
for every (resp. the given) path $\gamma $. In
this equality: $Y(\gamma )$ is $n n, n:=\dim(\pi ^{-1}(x)), x\in B$
nondegenerate depending only on $\gamma $ matrix; $Z(s;\gamma )$ is a
pseudo-unitary matrix of type $(p,q), q=n-p$, i.e.
\[
{\bigl(} Z(s;\gamma ){\bigr)}^{*} G_{p,q} Z(s;\gamma )
=G_{p,q}
:=
diag(\underbrace{1,\ldots,1}_{p\ times},
       \underbrace{-1,\ldots,-1}_{q\ times})
\qquad (4.10)
\]
 and $D(\gamma (s))$ is a fixed matrix such that
\[
{\bigl(}D(\gamma (s)){\bigr)}^{*}G^{T}(\gamma (s))D(\gamma
(s))=G_{p,q}.\qquad (4.11)
\]

{\bf Proof.} To prove the necessity we have to
solve the equation (4.6) with respect to $F$ when $G$ and $C$ are given. From
$G^{*}=G$ it follows $(G^{T})^{*}=G^{T}$which combined with the choice of
$\{e_{i}(s)\}$ and the independence of $p ($or/and $q)$ from $\gamma (s)$
leads to the existence of a satisfying (4.11) unitary matrix $D(\gamma (s))
[6]$.

Let $F(s;\gamma )=:F_{1}(s;\gamma ){\bigl(}D(\gamma (s)){\bigr)}^{-1}$. Then from (4.6) and (4.11), we get ${\bigl(}F^{-1}_{1}(s;\gamma ){\bigr)}^{*}G_{p,q}F^{-1}_{1}(s;\gamma )=C(\gamma )$. Putting here $F_{1}(s;\gamma )=:Y(\gamma )Z(s;\gamma )$, where $Y(\gamma )$ is an arbitrary nondegenerate (unitary) matrix for which $Y^{*}(\gamma )C(\gamma )Y(\gamma )=G_{p,q}$, we see that $Z(s;\gamma )$ satisfies (4.10). (The existence of $Y(\gamma )$ is a consequence of (4.6): from it and the inertial law of Jacobi-Sylvester it follows that $C(\gamma )$ has $p$ positive eigenvalues, hence the sought matrix $Y(\gamma )$ exists [6].) All this proves that $F(s;\gamma )$ has the form (4.9) under the condition that the considered $L$-transport along paths is consistent along every (resp.
some) path $\gamma $ with the given Hermitian bundle metric.

On the contrary, the sufficiency of the proposition is almost evident: if (4.9) is valid, then an elementary checking shows that (4.6) is true for $C(\gamma )={\bigl(}Y^{-1}(\gamma ){\bigr)}^{*}G_{p,q}Y^{-1}(\gamma )$ and according to proposition 4.2 the $L$-transport along paths and the Hermitian bundle metric are consistent along every (resp. some) path $\gamma .\blacksquare $

{\bf Proposition} ${\bf 4}{\bf .}{\bf 6}{\bf .} A$ necessary and sufficient
condition for the existence of a globally consistent with a given
$L$-transport along paths Hermitian bundle metrics is the existence of a
local basis in which the matrix $F(s;\gamma )$, defining this transport along
$\gamma :J\to B$ by (1.4){\bf ,} has the form (4.9) in which: $Y(\gamma )$
and $D(\gamma (s))$ are arbitrary nondegenerate $n n, n=\dim_{{\Bbb C}}(\pi
^{-1}(x)), x\in B$ matrices and $Z(s;\gamma )$ is a pseudo-unitary matrix of
type $(p,q), p+q=$n. Besides, if for a given $L$-transport along paths these
conditions are fulfilled, then in the above basis all globally consistent
with it Hermitian bundle metrics are described by the matrix
\[
G^{T}(\gamma (s))={\bigl(}(D(\gamma (s)))^{-1}
{\bigr)}^{*}G_{p,q}{\bigl(}D(\gamma (s)){\bigr)}^{-1},\qquad (4.12)
\]
which depends only on the point $\gamma (s)$, but not on the path $\gamma $.

{\bf Proof.} The proof of this proposition is an exact copy of the one of proposition 2.6 and it can be obtained from it with the following changes: metric$\mapsto $Hermitian metric; transposition sign $(T)\mapsto $Hermitian conjugation sign $(*); (2.8-10)\mapsto (4.8-10); G\mapsto  \mapsto G^{T}$and $G_{0}\mapsto G^{T}_{0}.\blacksquare $

\medskip
\medskip
 {\bf 5. REMARKS AND COMMENTS}

\medskip
(1) From the proof of proposition 2.1 it is clear that (2.2) is a matrix form of (1.3) in the basis $\{e_{i}(s)\}, s\in $J.

(2) In (2.8) the uncertainty in the choice of $F ($see (1.5)) is taken by $C ($see (2.7)), because of which $G(\gamma (s);\gamma )$ does not depend on the concrete choice of F.

(3) Said in another way, proposition 2.5 means that if there exist consistent with a given bundle metric $L$-transports along paths (along a given path), then they have in some basis a matrix (1.4) in which the matrix $F(s;\gamma )$ has the form (2.9).

(4) In proposition 3.1 the condition that the metric must be of class of smoothness $C^{1}$follows from the necessity for the equality (3.2) to have a sense and vice versa $(cf. [8,9])$.

(5) Proposition 3.4 follows also from proposition 3.3: it is sufficient to put $t=s_{0}$into (3.5) and  to denote $G(\gamma (t))$ with $C(\gamma )$.

(6) In a general form a necessary and sufficient condition for the existence of globally consistent with an $S$-transport metrics is given by proposition 2.6 in which $eq. (3.8)$ has to be taken into account, due to which the needed variant of proposition 2.6 can be formulated in terms of $\Gamma _{\gamma }$, but we are not going to do this here.

(7) In proposition 3.6 the condition for the independence of the signature of the metric from the point at which it is calculated is necessary for the existence of consistent with the metric $(S-)$transports along paths (see proposition 2.4).

(8) If $g_{i}(x)+g_{j}(x)\neq 0$ for every $i,j=1,\ldots  ,n$, of which type
are, in particular, the Euclidean metrics, then (3.13) may be written
equivalently as
\[
\Gamma _{\gamma }(s)=\Gamma _{1}(x)+\Gamma _{2}(x)G(x),
\qquad x=\gamma (s),\ s\in J,\qquad (5.1)
\]
 where $\Gamma _{2}=-\Gamma
^{T}_{2}$is arbitrary antisymmetric matrix and $\Gamma _{1}=\Gamma
^{T}_{1}$is arbitrary fixed symmetric solution of (3.6) with respect to
$\Gamma _{\gamma }$which under certain conditions (see below) admits the
representation (see [6], chapter $12, \sec. 13)$
\[
 \Gamma _{1}(\gamma (s))
=- \int\limits_{0}^{\infty}
{\bigl(}\exp(G(\gamma (s))t{\bigr)}\cdot
\frac{dG(\gamma(s)}{ds} \cdot
{\bigl(}\exp(G(\gamma (s))t){\bigr)}dt,\ s\in J.
 \qquad (5.2)
\]

Evidently, a necessary and sufficient condition for the existence of the
representation (5.1) is the existence of a symmetric solution $\Gamma
_{1}=\Gamma ^{T}_{1}$of the equation (3.6) with respect to $\Gamma _{\gamma
}$. The use of the described in the proof of proposition 3.6 method gives the
possibility to prove that such a solution exists iff for every pair $(i,j)$
for which $g_{i}(x)+g_{j}(x)=0$ the equality $K_{ij}(x)=0$ is satisfied
simultaneously. If (3.6) admits a symmetric solution $\Gamma _{1}$and there
exists at least one pair $(i,j)$ for which $g_{i}(x)+g_{j}(x)=K_{ij}(x)=0$,
then the integral in the right-hand side of (5.2) does not exist and $\Gamma
_{1}$admits the representation
\[
 \Gamma _{1}(x)=D(x)\Gamma _{0}(x)D^{T}(x), x\in \gamma (J),\qquad (5.3)
\]
where
\[
 \Gamma ^{T}_{0}:=\Gamma _{0}\qquad (5.4a)
\]
\[
{\bigl(}\Gamma _{0}(x){\bigr)}_{ij}
=K_{ij}(x)/(g_{i}(x)+g_{j}(x))
\qquad  for\
g_{i}(x)+g_{j}(x)\neq 0,\qquad (5.4b)
\]
and the remaining components of $\Gamma _{0}$, for every $(i,j)$ such that
$g_{i}(x)+g_{j}(x)=K_{ij}(x)=0$, if any, are arbitrary.

(9) The fact that all results for the consistency of real (symmetric,
Riemannian) metrics and linear transports along paths are true mutatis
mutandis also in the case of Hermitian metrics is not random. In fact, if we
denote by $h:x\mapsto h_{x}, x\in B$ an arbitrary Hermitian metric in the
complex vector bundle $(E,\pi ,B)$, then
\[
 g=Re(h)=\frac{1}{2}(h+h^{T}),\quad \omega =Im(h)=\frac{1}{2i} (h-h^{T}),
\ i=+\sqrt{-1},  \qquad (5.5)
\]
 where $h^{T}_{x}(u,v):=h_{x}(v,u), \pi
(u)=\pi (v)=x$, define, respectively, symmetric (Riemannian) and symplectic
metrics in $(E,\pi ,B)$. The definition of $h$ is equivalent to the
definition of $g$ or $\omega $, which is a corollary from
\[
 \omega _{x}(u,v)=g_{x}(u,{\bf J}v),\qquad (5.6)
\]
 where the complex
structure ${\bf J}$ on $(E,\pi ,B)$ is defined by ${\bf J}u=iu$. It is
important to note that
\[
 h=h\circ ({\bf J}\times {\bf J}),
\quad  g=g\circ ({\bf J}\times {\bf J}),
\quad \omega =\omega \circ ({\bf J}\times {\bf J}).\qquad (5.7)
\]
 Due to this the definition of an
arbitrary symmetric (complex) metric $g$ with the property
$g=g\circ ({\bf J}\times {\bf J})$ allows a Hermitian metric
$h=g+ig\circ (${\it id}$_{E} {\bf J})$ to be introduced. (In $[12] g$ itself
is called a Hermitian metric.) The existence of $g$ with the needed property
follows from the known fact that if $g_{0}$is a symmetric metric, then the
metric $g_{0}+g_{0}\circ ({\bf J}\times {\bf J})$ has the pointed property.

Namely the fact that the definition of a Hermitian metric $h$ is equivalent
to the definition of a symmetric metric $g$ with the property $g=g\circ ({\bf
J}\times {\bf J})$ is the reason that any result concerning symmetric metrics can
be formulated mutatis mutandis also for Hermitian metrics.

The above connections between Hermitian, Riemannian and symplectic metrics are not new and, for instance, can be found in the article "Hermitian Metric" in [13].

(10) The presented here material admits a generalization concerning arbitrary transports along paths in fibre bundles and bundle morphisms between them which will be a subject of other work.

\medskip
 {\bf ACKNOWLEDGEMENTS}

\medskip
The author expresses his gratitude to Prof. Vl. Aleksandrov (Institute of Mathematics of Bulgarian  Academy of Sciences) for constant interest in this work and stimulating discussions.

This research was partially supported by the Fund for Scientific Research of Bulgaria under contract Grant No. $F 103$.

\medskip
\medskip
 {\bf REFERENCES}

\medskip
1.  Greub W., S. Halperin, R. Vanstone, Connections, Curvature, and Cohomology, vol.1, vol.2, Academic Press, New York and London, $1972, 1973$.\par
2.  Iliev B.Z., Linear transports along paths in vector bundles. I. General theory, JINR Communication $E5-93-239$, Dubna, 1993.\par
3.  Dubrovin B.A., S.P. Novikov, A.T. Fomenko, Modern geometry, Nauka, Moscow, 1979 (In Russian).\par
4.  Synge J.L., Relativity: The general theory, North-Holland Publ. Co., Amsterdam, 1960.\par
5.  Iliev B.Z., Consistency between metrics and linear transports along curves, JINR Communication $E5-92-486$, Dubna, 1992.\par
6.  Bellman R., Introduction to matrix analysis, McGRAW-HILL book comp., New York-Toronto-London, 1960.\par
7.  Kobayashi S., K. Nomizu, Foundations of Differential Geometry, Vol. 1, Interscience Publishers, New York-London, 1963.\par
8.  Iliev B.Z., Parallel transports in tensor spaces generated by derivations of tensor algebras, JINR Communication $E5-93-1$, Dubna, 1993.\par
9.  Iliev B.Z., Linear transports along paths in vector bundles. II. Some  applications, JINR Communication $E5-93-260$, Dubna, 1993.\par
10.  Hicks N.J., Notes on Differential Geometry, D. Van Nostrand Comp., Inc., 1965.\par
11.  Karoubi M., $K$-theory. An Introduction, Springer-Verlag, Berlin-Heidelberg-New York, 1978.\par
12.  Kobayashi S., K. Nomizu, Foundations of differential geometry, vol. II, Interscience Publishers, New York-London, $1969, ch$. IX, \S4.\par
13.  Mathematical Encyclopedia, vol. 5, Moscow, Soviet encyclopedia, 1985 (In Russian).
\par
14.  Iliev B.Z., Transports along paths in fibre bundles. General theory, JINR Communication $E5-93-299$, Dubna, 1993.

\newpage

\vspace*{5ex}
\noindent
 Iliev B. Z. \\

\noindent
 Linear Transports along Paths in Vector Bundles \\
 IV. Consistency with Bundle Metrics\\[8ex]

 The problem for consistency between linear transports along paths and real
bundle metrics in real vector bundles is stated. Necessary and/or sufficient
conditions, as well as conditions for existence, for such consistency are
derived. All metrics (resp. transports) consistent with a given transport
(resp. metric) are explicitly obtained. The special case of generated by
derivations of tensor algebras linear transports of vectors is considered.
Analogous problems are investigated in complex vector bundles endowed with
Hermitian metrics.\\[5ex]

 The investigation has been performed at the Laboratory of Theoretical
Physics, JINR.

\end{document}